\documentclass[a4,12pt]{amsart}

\usepackage{amssymb}
\usepackage{epsf}
\usepackage{amsmath,epsf}
\usepackage[latin1]{inputenc}
\usepackage{amsfonts}

\numberwithin{equation}{section}
\newtheorem{theo}{Theorem}[section]
\newtheorem{prop}[theo]{Proposition}

\newtheorem{coro}[theo]{Corollary}

\newtheorem{rema}[theo]{Remark}
\newtheorem{exem}[theo]{Example}

\newenvironment{resume}{\footnotesize\quotation}

\newcommand{\pref}[1]{(\ref{#1})}

\def\cf{{\it cf. }}
\def\ie{{\it i.e. }}

\def\gx{>_{\rm lex}}

\def\nm{_{n,m}}

\def\C{{\mathbb C}}
\def\N{{\mathbb N}}

\def\Q{{\mathbb Q}}

\def\J{{\mathcal J}}

\def\S{{\mathcal S}}

\def\G{{\mathcal G}}

\def\B{{\mathcal B}}

\def\Qsym{{Qsym}}

\def\H{{\bf H}}
\def\SH{{\bf SH}}

\def\shuffle{{\,\raise
1pt\hbox{$\scriptscriptstyle\cup{\mskip-4mu}\cup$}\,}}

\title[Generalized Super-Coinvariant Polynomials]{Quasi-Invariant and Super-coinvariant Polynomials for the Generalized Symmetric Group}

\author{J.-C.~Aval}

\address[Jean-Christophe Aval]{LaBRI\\ Universit\'e Bordeaux 1\\ 351 cours
de 
 la Lib\'eration\\ 33405 Talence cedex\\ FRANCE}
\email{aval@labri.fr}

\date{\today\thanks{Research financed by EC's IHRP Programme, within the Research Training Network "Algebraic Combinatorics in Europe," grant HPRN-CT-2001-00272.}} 

\begin{document} 

\maketitle 

\begin{abstract} 
The aim of this work is to extend the study of super-coinvariant polynomials, introduced in \cite{a9,b1}, to the case of the generalized symmetric group $G_{n,m}$, defined as the wreath product $C_m\wr\S_n$ of the symmetric group by the cyclic group. We define a quasi-symmetrizing action of $G_{n,m}$ on $\Q[x_1,\dots,x_n]$, analogous to those defined in \cite{hivert} in the case of $\S_n$. The polynomials invariant under this action are called quasi-invariant, and we define super-coinvariant polynomials as polynomials orthogonal, with respect to a given scalar product, to the quasi-invariant polynomials with no constant term. Our main result is the description of a Gröbner basis for the ideal generated by quasi-invariant polynomials, from which we dedece that the dimension of the space of super-coinvariant polynomials is equal to $m^n\,C_n$ where $C_n$ is the $n$-th Catalan number. 
\end{abstract} 

\vskip 0.3cm
\begin{resume}
{\sc R\'esum\'e.} 
Le but de ce travail est d'étendre l'étude des polynômes super-coinvariants (définis dans \cite{a9}), au cas du groupe symétrique généralisé $G_{n,m}$, défini comme le produit en couronne $C_m\wr\S_n$ du groupe symétrique par le groupe cyclique. Nous définissons ici une action quasi-symétrisante de $G_{n,m}$ sur $\Q[x_1,\dots,x_n]$, analogue à celle définie dans \cite{hivert} dans le cas de $\S_n$. Les polynômes invariants sous cette action sont dits quasi-invariants, et les polynômes super-coinvariants sont les polynômes orthogonaux aux polynômes quasi-invariants sans terme constant (pour un certain produit scalaire). Notre résultat principal est l'obtention d'une base de Gröbner pour l'idéal engendré par les polynômes quasi-invariants. Nous en déduisons alors que la dimension de l'espace des polynômes super-coinvariants est $m^n\,C_n$ o\`u $C_n$ est le $n$-ième nombre de Catalan.
\end{resume}

\vskip 0.3cm
%%%%%%%%%%%%%%%%%%%%%%%%%%%%%%%%%%%%%%%%%%%%%%%%
\section{Introduction}

Let $X$ denote the alphabet in $n$ variables $(x_1,\dots,x_n)$ and $\C[X]$ denote the space of polynomials with complex coefficients in the alphabet $X$. Let $G\nm=C_m\wr\S_n$ denote the wreath product of the symmetric group $\S_n$ by the cyclic group $C_m$. This group is sometimes known as the {\sl generalized symmetric group} (\cf \cite{osima}). It may be seen as the group of $n\times n$ matrices in which each row and each column has exactly one non-zero entry (pseudo-permutation matrices), and such that the non-zero entries are $m$-th roots of unity. The order of $G\nm$ is $m^n\,n!$. When $m=1$, $G\nm$ reduces to the symmetric group $\S_n$, and when $m=2$, $G\nm$ is the hyperoctahedral group $B_n$, \ie the group of signed permutations, which is the Weyl group of type $B$ (see \cite{lus} for example for further details). The group $G\nm$ acts classically on $\C[X]$ by the rule
\begin{equation}\label{clas}
\forall g\in G\nm,\ \forall P\in\C[X],\ g.P(X)=P(X.{}^tg),
\end{equation}
where $g$ is the transpose of the matrix $g$ and $X$ is considered as a row vector.
Let
$$Inv\nm=\{P\in\C[X]\ /\ \forall g\in G\nm,\ g.P=P\}$$
denote the set of $G\nm$-invariant polynomials. Let us denote by $Inv\nm^+$ the set of such polynomials with no constant term. We consider the following scalar product on $\C[X]$:
\begin{equation}\label{scal}
\langle P,Q\rangle=P(\partial X)Q(X)\mid_{X=0}
\end{equation}
where $\partial X$ stands for $(\partial x_1,\dots,\partial x_n)$ and $X=0$ stands for $x_1=\cdots=x_n=0$. The space of $G\nm$-coinvariant polynomials is then defined by
\begin{eqnarray*}
Cov\nm&=&\{P\in\C[X]\ /\ \forall Q\in Inv\nm,\ Q(\partial X)P=0\}\\
&=&\langle Inv\nm^+\rangle^\perp\ \simeq\ \C[X]/\langle Inv\nm^+\rangle
\end{eqnarray*}
where $\langle S\rangle$ denotes the ideal generated by a subset $S$ of $\C[X]$.

A classical result of Chevalley \cite{che} states the following equality:
\begin{equation}
\dim Cov\nm=|G\nm|=m^n\, n!
\end{equation}
which reduces when $m=1$ to the theorem of Artin \cite{artin} that the dimension of the harmonic space $\H_n=Cov_{n,1}$ (\cf \cite{orbit}) is $n!$.

Our aim is to give an analogous result in the case of quasi-symmetrizing action. The ring $Qsym$ of quasi-symmetric functions was introduced by 
 Gessel~\cite{ges} as a source 
 of generating functions for $P$-partitions~\cite{stanley1} and  
 appears in more and more combinatorial 
   contexts~\cite{BMSW,stanley1,stanley2}. Malvenuto and 
 Reutenauer \cite{MR} proved a graded Hopf duality between $QSym$ and the  Solomon descent 
 algebras and Gelfand {\it et.~al.~}~\cite{NC} defined the graded Hopf 
 algebra $NC$ of non-commutative symmetric functions and identified it 
 with the Solomon descent algebra. 

In \cite{a9,b1}, Aval {\it et.~al.} investigated the space $\SH_n$ of super-coinvariant polynomials for the symmetric group, defined as the orthogonal (with respect to \pref{scal}) of the ideal generated by quasi-symmetric polynomials with no constant term, and proved that its dimension as a vector space equals the $n$-th Catalan number: 
\begin{equation}\label{catsh}
\dim\SH_n=C_n=\frac 1 {n+1} {2n \choose n}.
\end{equation}
Our main result is a generalization of the previous equation in the case of super-coinvariant polynomials for the group $G\nm$.

In Section 2, we define and study a ``quasi-symmetrizing'' action of $G\nm$ on $\C[X]$. We also introduce invariant polynomials under this action, which are called quasi-invariant, and polynomials orthogonal to quasi-invariant polynomials, which are called super-coinvariant. The Section 3 is devoted to the proof of our main result (Theorem \ref{main}), which gives the dimension of the space $SCov\nm$ of super-coinvariant polynomials for $G\nm$: we construct an explicit basis for $SCov\nm$ from which we deduce its Hilbert series.

%%%%%%%%%%%%%%%%%%%%%%%%%%%%%%%%%%%%%%%%%%%%%%%%%%%%%

\section{A quasi-symmetrizing action of $G\nm$}

We use vector notation for monomials. More precisely, for
$\nu=(\nu_1,\dots,\nu_n)\in\N^n$, we denote $X^\nu$ the monomial
\begin{equation}\label{mon}
x_1^{\nu_1} x_2^{\nu_2} \cdots x_n^{\nu_n}.
\end{equation}
For a polynomial $P\in\Q[X]$, we further denote $[X^\nu]\,P(X)$ as
the coefficient of the monomial $X^\nu$ in $P(X)$.

Our first task is to define a quasi-symmetrizing action of the group $G\nm$ on $\C[X]$, which reduces to the quasi-symmetrizing action of Hivert (\cf \cite{hivert}) in the case $n=1$. This is done as follows. Let $A\subset X$ be a subalphabet of $X$ with $l$ variables and $K=(k_1,\dots,k_l)$ be a vector of positive ($>0$) integers. If $B$ is a vector whose entries are distinct variables $x_i$ multiplied by roots of unity, the vector $(B)_<$ is obtained by ordering the elements in $B$ with respect to the variable order.
Now the quasi-symmetrizing action of $g\in G\nm$ is given by
\begin{equation}\label{quasa}
g\bullet A^K=w(g)^{c(K)} {(A.{}^t|g|)_<}^K
\end{equation}
where $w(g)$ is the weight of $g$, \ie the product of its non-zero entries, $|g|$ is the matrix obtained by taking the modules of the entries of $g$, and the oefficient $c(K)$ is defined as follows:
$$c(K)=\left\{
\begin{array}{ll}
0&{\rm if \ } \forall i,\ k_i\equiv 0\ [m]\\
1&{\rm if\ not.}
\end{array}\right.$$

\begin{exem}\rm If $m=3$ and $n=3$, and we denote by $j$ the complex number $j=e^{\frac{2i\pi}{3}}$, then for example
\begin{eqnarray*}
&&\left(
\begin{array}{ccc}
0&0&j\\
1&0&0\\
0&j&0
\end{array}
\right)\,\bullet\,(x_1^2\,x_2)\\
&=&(j^2)^1\left[\left(
\begin{array}{ccc}
0&0&1\\
1&0&0\\
0&1&0
\end{array}
\right)\,.\,(x_1,x_2)\right]_<^{(2,1)}\\
&=&j^2{(x_3,x_1)_<}^{(2,1)}\\
&=&j^2(x_1,x_3)^{(2,1)}\\
&=&j^2\,x_1^2\,x_3.
\end{eqnarray*}
\end{exem}

It is clear that this defines an action of the generalized symmetric group $G\nm$ on $\C[X]$, which reduces to Hivert's quasi-symmetrizing action (\cf \cite{hivert}, Proposition 3.4) in the case $m=1$.

Let us now study its invariant and coinvariant polynomials. We need to recall some definitions.

A {\sl composition} $\alpha =(\alpha_1,\alpha_2, \dots ,\alpha_k)$ of a
positive integer $d$ is an ordered list of positive integers ($>0$) whose
sum is $d$. For a {\sl vector} $\nu\in\N^n$, let $c(\nu)$ represent the composition
obtained by erasing zeros (if any) in $\nu$. A polynomial $P\in\Q[X]$ is said
to be {\sl quasi-symmetric} if and only if, for any $\nu$ and $\mu$ in
$\N^n$, we have
$$[X^\nu]P(X)=[X^\mu]P(X)$$
whenever $c(\nu)=c(\mu)$. The space of quasi-symmetric polynomials in $n$
variables is denoted
by $\Qsym_n$.

The polynomials invariant under the action \pref{quasa} of $G\nm$ are said to be {\it quasi-invariant} and the space of quasi-invariant polynomials is denoted by $QInv\nm$, \ie
$$P\in QInv\nm\ \Leftrightarrow\ \forall g\in G\nm,\ g\bullet P=P.$$
Let us recall (\cf \cite{hivert}, Proposition 3.15) that $QInv_{n,1}=QSym_n$. 
The following proposition gives a characterization of $QInv\nm$.
\begin{prop}\label{propu}
One has
$$P\in QInv_{n,m}\ \Leftrightarrow\ \exists Q\in QSym_n\ /\ P(X)=Q(X^m)$$
where $Q(X^m)=Q(x_1^m,\dots,x_n^m)$.
\end{prop}
\proof
Let $P$ be an element of $QInv\nm$. Let us denote by $\zeta$ the $m$-th root of unity $\zeta=e^{\frac{2i\pi} m}$ and by $g_j$ the element of $G\nm$ whose matrix is
$$\left(
\begin{array}{cccc}
\zeta&&0&\\
&1&&\\
0&&\ddots&\\
&&&1
\end{array}
\right)$$
with the $\zeta$ in place $j$.
Then we observe that the identities
$$\forall j=1,\dots,n,\ \frac 1 m (P+g_j\bullet P+g_j^2\bullet P+\cdots+g_j^{m-1}\bullet P)=P$$
imply that every exponents appearing in $P$ are multiples of $m$. Thus there exists a polynomial $Q\in\C[X]$ such that $P(X)=Q(X^m)$. 
To conclude, we note that $\S_n\subset G\nm$ implies that $P$ is quasi-symmetric, whence $Q$ is also quasi-symmetric.

The reverse implication is obvious.
\endproof

Let us now define {\it super-coinvariant} polynomials: 
\begin{eqnarray*}
SCov\nm&=&\{P\in\C[X]\ /\ \forall Q\in QInv\nm,\ Q(\partial X)P=0\}\\
&=&\langle QInv\nm^+\rangle^\perp\ \simeq\ \C[X]/\langle QInv\nm^+\rangle
\end{eqnarray*}
with the scalar product defined in \pref{scal}. This is the natural analogous to $Cov_n$ in the case of quasi-symmetrizing actions and $SCov\nm$ reduces to the space of super-harmonic polynomials $\SH_n$ (\cf \cite{b1}) when $m=1$. 

\begin{rema}\rm
It is clear that any polynomial invariant under \pref{quasa} is also invariant under \pref{clas}, \ie $Inv\nm\subset QInv\nm$. By taking the orthogonal, this implies that $SCov\nm\subset Cov\nm$. These observations somewhat justify the terminology.
\end{rema}

Our main result is the following theorem which is a generalization of equality \pref{catsh}.

\begin{theo}\label{main}
The dimension of the space $Scov\nm$ is given by
\begin{equation}\label{maineq}
\dim SCov\nm=m^n\,C_n=m^n\,\frac 1 {n+1} {2n \choose n}.
\end{equation}
\end{theo}

\begin{rema}\rm
In the case of the hyperoctahedral group $B_n=G_{n,2}$, C.-O. Chow \cite{chow} defined a class $BQSym(x_0,X)$ of quasi-symmetric functions of type $B$ in the alphabet $(x_0,X)$. His approach is quite different from ours. In particular, one has the equality:
$$BQSYm(x_0,X)=QSym(X)+QSym(x_0,X).$$
In the study of the coinvariant polynomials, it is not difficult to prove that the quotient $\C[x_0,X]/\langle BQSym^+\rangle$ is isomorphic to the quotient $\C[X]/\langle QSym^+\rangle$ studied in \cite{b1}. To see this, we observe that if $\G$ is the Gröbner basis of $\langle QSym^+\rangle$ constructed in \cite{b1} (see also the next section), then the set $\{x_0,\G\}$ is a Gröbner basis (any syzygy is reducible thanks to Buchberger's first criterion, \cf \cite{CLO}).
\end{rema}

The next section is devoted to give a proof of Theorem \ref{main} by constructing an explicit basis for the quotient $\C[X]/\langle QInv\nm^+\rangle$.

%%%%%%%%%%%%%%%%%%%%%%%%%%%%%%%%%%%%%%%%%%%%%%%%%%%%%

\section{Proof of the main theorem}

Our task is here to construct an explicit monomial basis for the quotient space $\C[X]/\langle QInv\nm^+\rangle$. Let us first recall (\cf \cite{b1}) the following bijection which associates to any vector $\nu\in\N^n$ a path $\pi(\nu)$ in the $\N\times\N$
plane with steps going north or east as follows. If
$\nu=(\nu_1,\dots,\nu_n)$, the path $\pi(\nu)$ is
$$(0,0)\rightarrow(\nu_1,0)\rightarrow(\nu_1,1)\rightarrow(\nu_1+\nu_2,1)
\rightarrow(\nu_1+\nu_2,2)\rightarrow\cdots\ \ \ \ \ \ \ \ \ \ \ \ \ \ \ \ \ \
\ \ \ $$
\vskip -0.6cm
$$\ \ \ \ \ \ \ \ \ \ \ \ \ \ \ \ \ \ \ \ \ \ \ \ \ \ \ \ \ \ \ \ \ \ \ \ \
\rightarrow(\nu_1+\cdots+\nu_n,n-1)\rightarrow(\nu_1+\cdots+\nu_n,n).$$
For example the path associated to $\nu=(2,1,0,3,0,1)$ is

\vskip 0.2cm
\centerline{\epsffile{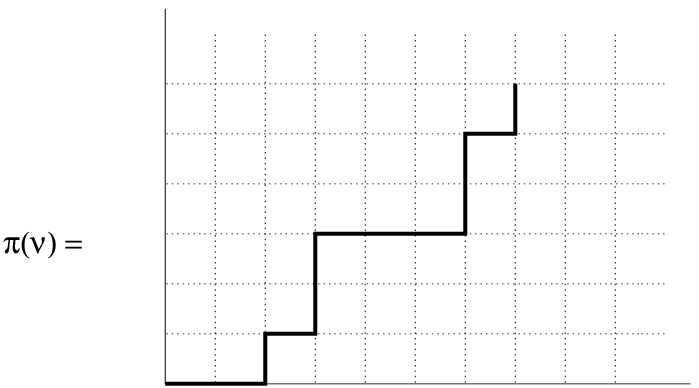}}

\vskip 0.2cm

We distinguish two kinds of paths, thus two kinds of vectors, with respect
to their ``behavior'' regarding the diagonal $y=x$.
If the path remains above the diagonal, we call it a {\sl Dyck path}, and
say that the corresponding vector is {\sl Dyck}. If not, we say that the
path (or equivalently the associated vector) is {\sl transdiagonal}. For
example $\eta=(0,0,1,2,0,1)$ is Dyck and $\varepsilon=(0,3,1,1,0,2)$ is
transdiagonal.

\vskip 0.2cm
\centerline{\epsffile{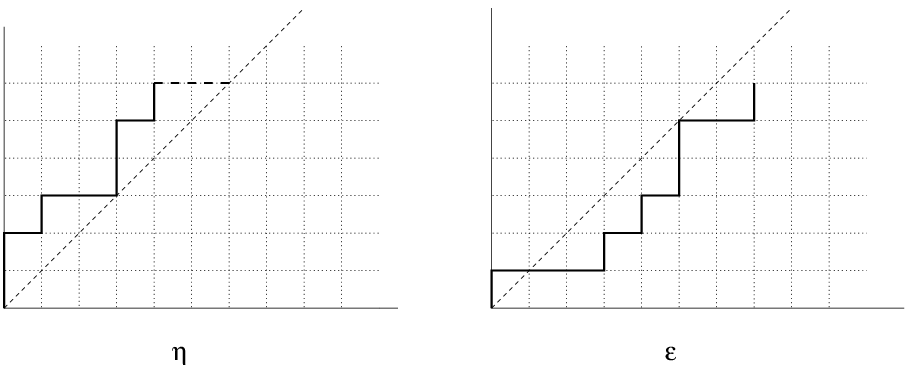}}

\vskip 0.2cm

We then have the following result which generalizes Theorem 4.1 of \cite{b1} and which clearly implies the Theorem \ref{main}.
\begin{theo}\label{main2}
The set of monomials
$$\B\nm=\{(X_n)^{m \,\eta+\alpha} /\ \pi(\eta)\ {\it is\ a\ Dyck\ path},\ 0\le\alpha_i<m\}$$
is a basis for the quotient $\C[X_n]/\langle QInv_{n,m}^+\rangle$.
\end{theo}

To prove this result, the goal is here to construct a Gröbner basis for the ideal $\J\nm=\langle QInv_{n,m}^+\rangle$. We shall use results of \cite{a9,b1}.

Recall that the {\it lexicographic order} on monomials
is
\begin{equation}
X^\nu\gx X^\mu\qquad{\rm iff}\qquad \nu\gx \mu,
\end{equation}
if and only if the first non-zero part of the
vector $\nu-\mu$ is positive.

For any subset $\S$ of $\Q[X]$ and for any positive integer $m$, let us introduce $\S^m=\{P(X^m)\ ,\ P\in\S\}$. If we denote by $G(I)$ the unique reduced monic Gröbner basis (\cf \cite{CLO}) of an ideal $I$, then the simple but crucial fact in our context is the following.

\begin{prop}\label{ppp}
With the previous notations,
\begin{equation}
G(\langle\S^m\rangle)=G(\langle\S\rangle)^m.
\end{equation}
\end{prop} 
\proof
This is a direct consequence of Buchberger's criterion. Indeed, if for every pair $g,\,g'$ in $G(\langle S\rangle)$, the syzygy
$$S(g,g')$$
reduces to zero, then the syzygy
$$S(g(X^m),g'(X^m))$$
also reduces to zero in $G(\langle S^m\rangle)$ by exactly the same computation.
\endproof

Let us recall that in \cite{a9} is constructed a family $\G$ of polynomials $G_\varepsilon$ indexed by transdiagonal vectors $\varepsilon$. This family is constructed by using recursive relations of the fundamental quasi-symmetric functions and one of its property (\cf \cite{a9}) says that the leading monomial of $G_\varepsilon$ is: $LM(G_\varepsilon)=X^\varepsilon$. 
Since $\G$ is a Gröbner basis of $\J_{n,1}$, the following result is a consequence of Propositions \ref{propu} and \ref{ppp}. 
\begin{prop}
The set $\G^m$ is a Gröbner basis of the ideal $\J\nm$.
\end{prop}

To conclude the proof of Theorem \ref{main2}, it is sufficient to observe that the monomials not divisible by a leading monomial of an element of $\G^m$, \ie by a $X^{m\varepsilon}$ for $\varepsilon$ transdiagonal, are precisely the monomials appearing in the set $\B\nm$.

As a corollary of Theorem \ref{main2}, one gets an explicit formula for the Hilbert series of $SCov\nm$. For $k\in\N$, let $SCov\nm^{(k)}$ denote the projection
\begin{equation}
SCov\nm^{(k)}=SCov\nm\,\cap\,\Q^{(k)}[X]
\end{equation}
where $\Q^{(k)}[X]$ is the vector space of homogeneous polynomials of degree
$k$ together with zero. 

Let us denote by $F\nm(t)$ the Hilbert series of $SCov\nm$, \ie
\begin{equation}
F\nm(t)=\sum_{k\ge 0} \dim SCov\nm^{(k)}\, t^k.
\end{equation}

Let us recall that in \cite{b1} is given an explicit formula for $F_{n,1}$:
\begin{equation}\label{ff}
F_{n,1}(t)=F_n(t)=\sum_{k=0}^{n-1}\frac{n-k}{n+k}{n+k
\choose k} t^k
\end{equation}
using the number of Dyck paths with a given number of factors (\cf \cite{Kreweras}).

The Theorem \ref{main2} then implies the

\begin{coro}
With the notations of \pref{ff}, the Hilbert series of $SCov\nm$ is given by
$$F\nm(t)=\frac{1-t^m}{1-t}F_n(t^m)$$
from which one deduces the close formula
$$\sum_n F_{n,m}(t)\,x^n=\frac{(1\!-\!t)-\sqrt{(1\!-\!t)(1-t-4t^mx(1-t^m))}-2x(1\!-\!t^m)}{(1-t)(2t^m-1)-x(1-t^m)}\raise 2pt\hbox{.}$$
\end{coro}

%%%%%%%%%%%%%%%%%%%%%%%%%%%%%%%%%%%%%%%%%%%%%%%%%%%%%

\end{document}